\documentclass[11pt]{amsart}
\usepackage{amssymb,latexsym,amsmath,amsfonts}
\usepackage[mathscr]{eucal}
\numberwithin{equation}{section}
\theoremstyle{definition}
\newtheorem{definition}{Definition}[section]
\newtheorem{example}[definition]{Example}
\newtheorem{custom}[definition]{}
\theoremstyle{remark}
\newtheorem{remark}[definition]{Remark}
\theoremstyle{plain}
\newtheorem{theorem}[definition]{Theorem}
\newtheorem{lemma}[definition]{Lemma}

\newtheorem{result}[definition]{Result}

\newcommand{\zt}{\zeta}
\newcommand{\zbar}{\overline{z}}

\newcommand{\lam}{\lambda}
\newcommand{\exep}{\mathfrak{S}}

\newcommand{\banal}{\mathcal{A}}
\newcommand{\id}{\mathbb{I}}

\newcommand{\bdy}{\partial D}

\newcommand{\OM}{\Omega}

\newcommand{\hol}{\mathcal{O}}

\newcommand\hyper[2]{\left|\frac{{#1}-{#2}}{1-\overline{{#2}}{#1}}\right|}
\newcommand\ilnhyper[2]{|({#1}-{#2})(1-\overline{{#2}}{#1})^{-1}|}

\newcommand\blah[2]{\left(\frac{{#1}-{#2}}{1-\overline{{#2}}{#1}}\right)}
\newcommand{\mobi}{\mathcal{M}}
\newcommand{\mobd}{{\sf dist}_{\mathcal{M}}}
\newcommand{\impl}{\Longrightarrow}
\newcommand{\mapp}{\longrightarrow}

\newcommand{\Fd}{F_d}

\newcommand{\opti}{\mathscr{F}^{z,w}}
\newcommand{\shrp}{\mathfrak{G}^{A,d}}
\newcommand{\Shrp}{\mathfrak{G}^{(d)}}

\newcommand{\tr}{{\sf tr}}

\newcommand{\cplx}{\mathbb{C}}

\begin{document}

\title[Schwarz lemmas for the spectral unit ball]{Sharp Schwarz-type lemmas for the \\ 
spectral unit ball}
\author{Gautam Bharali}
\address{Department of Mathematics, Indian Institute of Science, Bangalore -- 560 012}
\email{bharali@math.iisc.ernet.in}

\begin{abstract} 
We provide generalisations of two Schwarz-type lemmas --- the first a result of Globevnik and the 
other due to Ransford and White ---  for holomorphic mappings into the spectral unit ball. The first
concerns mappings of the unit disc in $\cplx$ into the spectral unit ball, while the second concerns
self-mappings. The aforementioned results apply to holomorphic mappings that map the origin to 
the origin. We extend these results to general holomorphic mappings into the spectral unit ball.
We also show that our results are sharp.
\end{abstract}
\maketitle

\section{Introduction and Statement of Results}\label{S:intro}

The {\bf spectral unit ball} (denoted by $\OM_n$) is defined as
\[
\OM_n \ := \ \{W\in M_n(\cplx):r(W)<1\},
\]
where $r(W)$ denotes the spectral radius of the $n\times n$ matrix $W$. In this paper, $D$ will 
denote the unit disc in $\cplx$ and, if $\OM$ and $G$ are complex domains, $\hol(\OM;G)$ will denote
the class of holomorphic mappings of $\OM$ into $G$. We present two Schwarz-type lemmas for the 
spectral unit ball. These lemmas are inspired by the renewed interest in the function
theory on the spectral unit ball --- the reader is referred to \cite{baribeauRansford:nlspm00},
\cite{rostand:asub03}, \cite{aglerYoung:2b2sNPp04} and \cite{edgrnZwnk:gsp05}, to name
just a few recent papers. This interest stems, to a large extent, from recent work on 
the spectral version of the Nevanlinna-Pick
interpolation problem. We shall not address this problem directly in this paper; although 
Theorem~\ref{T:disc} below might have some bearing on the {\em two-point} interpolation problem.
We begin by considering the following result by Globevnik \cite{globevnik:Slsr74}, which is perhaps
the earliest  Schwarz-type lemma for the spectral unit ball:
\begin{equation}\label{E:globevnik}
F\in\hol(D;\OM_n) \ \text{and} \ F(\zt_1)=0 \ \impl \ 
	r(F(\zt_2))\leq\hyper{\zt_1}{\zt_2}=:\mobi(\zt_1,\zt_2).
\end{equation}
One would like to generalise this result to the case when $F(\zt_j)$ is not necessarily $0$,
$j=1,2$. 
\smallskip

The second Schwarz-type lemma that motivates this paper --- this time a result on
{\em self-mappings} of $\OM_n$ --- is the following result of Ransford and White 
\cite{ransfordWhite:hsmsub91}:
\begin{equation}\label{E:rnsfrdWht}
G\in\hol(\OM_n;\OM_n) \ \text{and} \ G(0)=0 \ \impl \ r(G(X))\leq r(X) \;\; \forall X\in\OM_n.
\end{equation}
Incidentally, we refer to the above results as ``Schwarz-type lemmas'' because they relate 
the growth of a holomorphic mapping to the growth of its argument(s). One would like 
to generalise \eqref{E:rnsfrdWht} in the similar manner that the Schwarz-Pick lemma
generalises the Schwarz lemma for $D$ --- i.e. by formulating an inequality that is valid without
assuming that the holomorphic mapping in question has a fixed point.
\smallskip

What is the key idea needed to generalise \eqref{E:globevnik} and \eqref{E:rnsfrdWht} 
in appropriate ways~? The following example shows how the conclusion of \eqref{E:globevnik}
can fail, even in the simple situation where $r(F(\zt_1))=0$, if $F(\zt_1)\neq 0$. However,
it also suggests a way forward.
\smallskip

\begin{example}\label{Ex:keyex}
For $n\geq 3$ and $d=2,\dots,n-1$, define the holomorphic map $\Fd:D\mapp\OM_n$ by
\[
\Fd(\zt) \ := \ \begin{bmatrix}
                \ 0  & {} & {} & \zt \ &\vline & \ {} \ \\
                \ 1  & 0  & {} & 0 \ &\vline & {} \ \\
                \ {} & \ddots  & \ddots & \vdots \ &\vline & \text{\LARGE{0}} \ \\  
                \ {} & {} & 1 & 0 \ &\vline & {} \ \\ \hline
                \ {} & {} & {} & {} \ &\vline & {} \ \\
                \ {} & {} & \text{\LARGE{0}} & {} \ &\vline & \zt\mathbb{I}_{n-d} \
                \end{bmatrix}_{n\times n}, \qquad \zt\in D,
\]
where $\mathbb{I}_{n-d}$ denotes the identity matrix of dimension $n-d$ for $1<d<n$. 
One easily computes that $r(\Fd(\zt))=|\zt|^{1/d}$. Hence
\begin{align}
r(\Fd(\zt))^d \ &= \ |\zt|\ =\ \mobi(0,\zt)\quad \forall\zt\in D, \label{E:nilpot} \\
\text{but, for each $q<d$,} \ r(\Fd(\zt))^q \ &> \ \mobi(0,\zt)\quad \forall\zt\neq 0. \notag
\end{align}
In particular, $r(\Fd(\zt))>\mobi(0,\zt)\quad \forall\zt\neq 0$, in contrast with 
\eqref{E:globevnik}.\qed
\end{example}
\smallskip

The above example is rather suggestive when one notices that the exponent occurring in
the left-hand side of \eqref{E:nilpot} is the degree of the minimal polynomial of $\Fd(0)$. While
$r(F_d(\zt_1))=0$ for each $d=2,\dots,n-1$ (take $\zt_1=0$ in this discussion), what differs in
each case is the degree of the minimal polynomial of $F_d(\zt_1)$. Presumably, this information 
should be encoded in any generalisation of \eqref{E:globevnik}. This idea is the key to establishing
a result with the following features:
\begin{itemize}
\item It has a Schwarz-type structure:  i.e. we get an expression in $F(\zt_1)$ and 
$F(\zt_2)$ --- call it $\mathcal{E}(F(\zt_1),F(\zt_2))$ --- such that
\[
F\in\hol(D;\OM_n) \impl \mathcal{E}(F(\zt_1),F(\zt_2)) \ \leq \ \mobi(\zt_1,\zt_2).
\]
\item Globevnik's result is recovered when we set $F(\zt_1)=0$ in the above inequality.
\item The above Schwarz-type inequality is {\em sharp} in the sense that this inequality
is the best one can achieve. In more precise terms: given $z\neq w\in D$, we can find a 
$\opti\in\hol(D;\OM_n)$ such that $\mathcal{E}(\opti(z),\opti(w))=\mobi(z,w)$.
\end{itemize}
Before stating this result, let us, for any compact subset $K\varsubsetneq D$ and 
$\zt\in D$, define
\[
\mobd(\zt;K) \ := \ \min_{z\in K}\hyper{\zt}{z}.
\]
Our first result is as follows.
\smallskip

\begin{theorem}\label{T:disc} Let $F\in\hol(D;\OM_n)$, $n\geq 2$, and let $\zt_1,\zt_2\in D$. Let
$W_j=F(\zt_j), \ j=1,2$, and define
\[
d_j \ := \ \text{the degree of the minimal polynomial of $W_j$},
\]
for $j=1,2$. Then
\begin{equation}\label{E:SchwarzIneq}
\max\left\{\max_{\mu\in\sigma(W_2)}[\mobd(\mu;\sigma(W_1))^{d_1}], \ \max_{\lambda\in\sigma(W_1)}
[\mobd(\lambda;\sigma(W_2))^{d_2}]\right\} \ \leq \ \hyper{\zt_1}{\zt_2}.
\end{equation}
\smallskip

Furthermore, \eqref{E:SchwarzIneq} is sharp in the sense that given any two points $z,w\in D$,
there exists a mapping $\opti\in\hol(D;\OM_n)$ such that
\begin{multline}
\max\left\{\max_{\mu\in\sigma(\opti(w))}[\mobd(\mu;\sigma(\opti(z)))^{d(z)}], \ 
\max_{\lambda\in\sigma(\opti(z))}[\mobd(\lambda;\sigma(\opti(w)))^{d(w)}]\right\} \\
= \ \hyper{z}{w},\notag
\end{multline}
where $d(z)$ (resp.~$d(w)$) is the degree of the minimal polynomial of $\opti(z)$ 
(resp.~$\opti(w)$).
\end{theorem}

\begin{remark} Note that Globevnik's result is recovered when we set $F(\zt_1)=0$ in 
the above theorem. This is because if $W_1=0$, then, in the notation of Theorem~\ref{T:disc},
$d_1=1$, and
\begin{align}
\max_{\lambda\in\sigma(W_1)}[\mobd(\lambda;\sigma(W_2))^{d_2}] \ 
	&= \ \min_{\mu\in\sigma(W_2)}|\mu|^{d_2}\notag \\ 
&\leq \ \max_{\mu\in\sigma(W_2)}[\mobd(\mu;\sigma(W_1))^{d_1}] \ 
	= \ r(W_2).\notag
\end{align}
Hence \eqref{E:SchwarzIneq}, in this case, is identical to the conclusion of Globevnik's result.
\end{remark}
\smallskip

The proof of Theorem~\ref{T:disc} is presented in Section~\ref{S:proofDisc}. The consideration
of a pertinent minimal polynomial turns out to equally relevant to our next result. Essentially,
the proofs of both results exploit the minimal polynomial of a key matrix lying in the ranges of $F$,
respectively $G$, to transform the maps $F$ and $G$ to maps to which the results \eqref{E:globevnik} 
and \eqref{E:rnsfrdWht}, respectively, are applicable. One is led to do this when examining the 
basic example of a mapping 
$G\in\hol(\OM_n;\OM_n)$ where $G(0)\neq 0$ and $G$ fails to satisfy the inequality in 
\eqref{E:rnsfrdWht}, {\em even though} $r(G(0))=0$. Since we do not wish to prolong
this already protracted introduction, we refer the reader to the counterexample immediately
following Theorem 2 in the Ransford-White paper \cite{ransfordWhite:hsmsub91} (or to the 
end of Section~\ref{S:proofSpec} of this paper). The idea hinted at above leads to a new result that:
\begin{itemize}
\item Just like \eqref{E:rnsfrdWht}, provides a bound on the growth of the spectral radius of $G(X)$
in terms of $r(X)$, and specialises precisely to \eqref{E:rnsfrdWht} when we set $G(0)$;
\item Is sharp in a manner analogous to our discussion of ``sharpness'' of our previous result.
\end{itemize}
More precisely, we have the following:
\smallskip

\begin{theorem}\label{T:spec} Let $G\in\hol(\OM_n;\OM_n), \ n\geq 2$, and define $d_G:=$ the 
degree of the minimal polynomial of $G(0)$. Then:
\begin{equation}\label{E:growthBound}
r(G(X)) \ \leq \ \frac{r(X)^{1/d_G}+r(G(0))}{1+r(G(0))r(X)^{1/d_G}} \quad\forall X\in\OM_n.
\end{equation}
\smallskip

Furthermore, the inequality \eqref{E:growthBound} is sharp in the sense that there exists a 
non-empty set $\exep_n\subset\OM_n$ such that given any $A\in\exep_n$ and $d=1,\dots,n$, we can
find a $\shrp\in\hol(\OM_n;\OM_n)$ such that
\begin{align}
d_{\shrp} \ &= \ d, \;\; \text{and}\notag\\
r(\shrp(A)) \ &= \ \frac{r(A)^{1/d}+r(\shrp(0))}{1+r(\shrp(0))r(A)^{1/d}}. \label{E:sharp}
\end{align}
\end{theorem}

\begin{remark} It is quite obvious why the Ransford-White bound is recovered when we set
$G(0)=0$ in the above theorem. When $G(0)=0$, then $d_G=1$ and $r(G(0))=0$, whence
\eqref{E:growthBound} is identical to the conclusion of \eqref{E:rnsfrdWht}.
\end{remark}
\medskip

\section{The Proof of Theorem~\ref{T:disc}}\label{S:proofDisc}

The proofs in this section depend crucially on a theorem by Vesentini. The result is as follows:
  
\begin{result}[Vesentini, \cite{vesentini:spr68}]\label{R:vesentini}
Let $\banal$ be a complex, unital Banach
algebra and let $r(x)$ denote the spectral radius of any element $x\in\banal$.
Let $f\in\hol(D;\banal)$. The the function $\zt\longmapsto r(f(\zt))$ is
subharmonic on $D$.
\end{result}
\smallskip

The following result is the key lemma of this section. The proof of Theorem~\ref{T:disc} is reduced to
a simple application of this lemma.
\smallskip

\begin{lemma}\label{L:key} Let $F\in\hol(D;\OM_n)$ and let $\lam_1,\dots,\lam_s$ be the distinct
eigenvalues of $F(0)$. Define $m(j):=$the multiplicity of the factor $(\lam-\lam_j)$ in the
minimal polynomial of $F(0)$. Define the Blaschke product
\[
B(\zt) \ := \ \prod_{j=1}^s\blah{\zt}{\lam_j}^{m(j)}, \quad\zt\in D.
\]
Then $|B(\lam)|\leq|\zt| \ \forall\lam\in\sigma(F(\zt))$.
\end{lemma}
\begin{proof}
The Blaschke product $B$ induces a matrix function $\widetilde{B}$ on $\OM_n$: for any matrix 
$A\in\OM_n$, we set
\[
\widetilde{B}(A) \ := \ \prod_{j=1}^s(\id-\overline{\lam_j} A)^{-m(j)}(A-\lam_j\id)^{m(j)},
\]
which is well-defined on $\OM_n$ because whenever $\lam_j\neq 0$,
\[
(\id-\overline{\lam}_j A) \ = \ \overline{\lam_j}(\id/\overline{\lam}_j-A)\in GL(n,\cplx).
\]
Furthermore, since $\zt\longmapsto(\zt-\lam_j)/(1-\overline{\lam_j}\zt)$ has a power-series
expansion that converges uniformly on  compact subsets of $D$, it follows from standard arguments
that
\begin{equation}\label{E:specB}
\sigma(\widetilde{B}(A)) \ = \ \{B(\lam):\lam\in\sigma(A)\}\quad\text{for any $A\in\OM_n$.}
\end{equation}
By the definition of the minimal polynomial, $\widetilde{B}\circ F(0)=0$. At this point, we could
apply Globevnik's lemma --- i.e. \eqref{E:globevnik} above --- to complete the proof. The actual
argument, however, is very elementary, and we provide it here. Since $\widetilde{B}\circ F(0)=0$, there
exists a holomorphic map $\Phi\in\hol(D;M_n(\cplx))$ such that $\widetilde{B}\circ F(\zt)=\zt\Phi(\zt)$.
Note that
\begin{equation}\label{E:specRelation}
\sigma(\widetilde{B}\circ F(\zt)) \ = \ \sigma(\zt\Phi(\zt)) \ = \ 
\zt\sigma(\Phi(\zt))\quad\forall\zt\in D.
\end{equation}
Since $\sigma(\widetilde{B}\circ F(\zt))\subset D$, the above equations give us:
\begin{equation}\label{E:circBound}
r(\Phi(\zt)) \ < \ 1/R \quad\forall\zt:|\zt|=R, \ R\in(0,1).
\end{equation}
Taking $\banal=M_n(\cplx)$ in Vesentini's theorem, we see that
$\zt\longmapsto r(\Phi(\zt))$ is subharmonic on the unit disc. Applying the Maximum Principle to
\eqref{E:circBound} and taking limits as $R\longrightarrow 1^-$, we get
\begin{equation}\label{E:oneBound}
r(\Phi(\zt)) \ \leq \ 1\quad\forall\zt\in D.
\end{equation}
In view of \eqref{E:specB}, \eqref{E:specRelation} and \eqref{E:oneBound}, we get
\[
|B(\lam)| \ \leq \ |\zt|r(\Phi(\zt)) \ \leq \ |\zt|\quad\forall\lambda\in\sigma(F(\zt)).
\]
\end{proof}
\smallskip

We are now in a position to provide
\smallskip

\begin{custom}\begin{proof}[{\bf The proof of Theorem~\ref{T:disc}:}] Define the 
disc automorphisms
\[
M_j(\zt) \ := \ \frac{\zt-\zt_j}{1-\overline{\zt_j}\zt}, \quad j=1,2,
\]
and write $\Phi_j=F\circ M_j^{-1}, \ j=1,2$. Note that $\Phi_1(0)=W_1$. Let 
$\lam_1,\dots,\lam_r$ be the distinct eigenvalues of $W_1$ and define
$m_1(j):=$the multiplicity of the factor $(\lam-\lam_j)$ in the minimal polynomial
of $W_1$. Define
\[
B_1(\zt) \ := \ \prod_{j=1}^r\blah{\zt}{\lam_j}^{m_1(j)}, \quad\zt\in D.
\]
Applying Lemma~\ref{L:key}, we get
\begin{align}
\hyper{\zt_1}{\zt_2} \ = \ |M_1(\zt_2)| \ &\geq \ 
\prod_{j=1}^r\hyper{\mu}{\lam_j}^{m_1(j)} \label{E:1stpart} \\
&\geq \ \mobd(\mu;\sigma(W_1))^{d_1} \quad\forall\mu\in\sigma(\Phi_1(M_1(\zt_2)))
							=\sigma(W_2).\notag
\end{align}
Now, swapping the roles of $\zt_1$ and $\zt_2$ and applying the same argument to
\[
B_2(\zt) \ := \ \prod_{j=1}^s\blah{\zt}{\mu_j}^{m_2(j)}, \quad\zt\in D,
\]
where $\mu_1,\dots,\mu_s$ are the distinct eigenvalues of $W_2$ and $m_2(j):=$the multiplicity
of the factor $(\lam-\mu_j)$ in the minimal polynomial of $W_2$, we get
\begin{equation}\label{E:2ndpart}
\hyper{\zt_1}{\zt_2} \ \geq \
\mobd(\lam;\sigma(W_2))^{d_2} \quad\forall\lam\in\sigma(W_1).
\end{equation}
Combining \eqref{E:1stpart} and \eqref{E:2ndpart}, we get
\[
\max\left\{\max_{\mu\in\sigma(W_2)}[\mobd(\mu;\sigma(W_1))^{d_1}], \ \max_{\lambda\in\sigma(W_1)}
[\mobd(\lambda;\sigma(W_2))^{d_2}]\right\} \ \leq \ \hyper{\zt_1}{\zt_2}.
\]

In order to prove the sharpness of \eqref{E:SchwarzIneq}, fix an $n\geq 2$, and choose any 
$z,w\in D$. Next, define $M(\zt):=(\zt-z)(1-\zbar\zt)^{-1}$. Pick any $d=1,\dots,n$, and define
\[
N_d(\zt) \ := \ \begin{cases}
		\quad [M(\zt)],        & \text{if $d=1$}, \\
		\ \begin{bmatrix}
			\ 0  & {} & {} & M(\zt) \ \\
			\ 1  & 0  & {} & 0 \ \\
			\ {} & \ddots & \ddots & \vdots \ \\
			\ {} & {} & 1 & 0 \
			\end{bmatrix}_{d\times d}, & \text{if $d\geq 2$},
		\end{cases}
\]
and, for the chosen $d$, define $\opti$ by the following block-diagonal matrix
\[
\opti(\zt) \ := \begin{bmatrix}
			\ N_d(\zt) & {} \ \\
			\ {} & M(\zt)\mathbb{I}_{n-d} \ 
			\end{bmatrix} \quad\forall\zt\in D.
\]
Note that
\begin{itemize}
\item $\opti(z)$ is nilpotent of degree $d$, whence $d(z)=d$; and
\item Since $|M(w)|^{1/d}>|M(w)|$,
\[
\max_{\mu\in\sigma(\opti(w))}[\mobd(\mu;\sigma(\opti(z)))^{d(z)}] \ = \ |M(w)| \ = 
\ \hyper{z}{w}.
\]
\end{itemize}
A similar argument yields 
\[
\max_{\lambda\in\sigma(\opti(z))}[\mobd(\lambda;\sigma(\opti(w)))^{d(w)}] \ = \ |M(w)|^{d+1} \
= \ \hyper{z}{w}^{d+1}.
\]
Hence, we have the equality
\begin{multline}
\max\left\{\max_{\mu\in\sigma(\opti(w))}[\mobd(\mu;\sigma(\opti(z)))^{d(z)}], \
\max_{\lambda\in\sigma(\opti(z))}[\mobd(\lambda;\sigma(\opti(w)))^{d(w)}]\right\} \\
= \ \hyper{z}{w}.\notag
\end{multline}
$\opti$ is therefore the desired map that establishes the sharpness of \eqref{E:SchwarzIneq}.		
\end{proof}
\end{custom}
\medskip

\section{The Proof of Theorem~\ref{T:spec}}\label{S:proofSpec}

In order to prove Theorem~\ref{T:spec}, we shall need the following elementary
\smallskip

\begin{lemma}\label{L:mobiusTrans} Given a M{\"o}bius transformation $T(z):=(az+b)/(cz+d)$,
if $T(\bdy)\Subset\cplx$, then $T(\bdy)$ is a circle with
\[
{\rm centre}(T(\bdy)) \ = \ \frac{b\overline{d}-a\overline{c}}{|d|^2-|c|^2}, \qquad
{\rm radius}(T(\bdy)) \ = \ \frac{|ad-bc|}{||d|^2-|c|^2|}.
\]
\end{lemma}
\smallskip

We are now in a position to present
\smallskip

\begin{custom}\begin{proof}[{\bf The proof of Theorem~\ref{T:spec}:}] Let $G\in\hol(\OM_n;\OM_n)$ 
and let $\lam_1,\dots,\lam_s$ be the distinct eigenvalues of $G(0)$. Define $m(j):=$the multiplicity 
of the factor $(\lam-\lam_j)$ in the minimal polynomial of $G(0)$. Define the Blaschke product
\[
B_G(\zt) \ := \ \prod_{j=1}^s\blah{\zt}{\lam_j}^{m(j)}, \quad\zt\in D.
\]
$B_G$ induces the following matrix function which, by a mild abuse of notation, we shall also denote 
as $B_G$
\[
B_G(Y) \ := \ \prod_{j=1}^s(\id-\overline{\lam_j} Y)^{-m(j)}(Y-\lam_j\id)^{m(j)} \quad\forall Y\in\OM_n,
\]
which is well-defined on $\OM_n$ precisely as explained in the proof of Lemma~\ref{L:key}. Once again,
owing to the analyticity of $B_G$ on $D$,
\[
\sigma(B_G(Y)) \ = \ \{B_G(\lam):\lam\in\sigma(Y)\}\quad\forall Y\in\OM_n,
\]
whence $B_G:\OM_n\mapp\OM_n$. Therefore, if we define
\[
H(X) \ := \ B_G\circ G(X) \quad\forall X\in\OM_n,
\]
then $H\in\hol(\OM_n;\OM_n)$ and, by construction, $H(0)=0$. By the Ransford-White result,
$r(H(X))\leq r(X)$, or, more precisely
\[
\max_{\mu\in G(X)}\left\{\prod_{j=1}^s\hyper{\mu}{\lam_j}^{m(j)}\right\} \ \leq \ 
r(X) \quad\forall X\in\OM_n.
\]
In particular:
\[
\max_{\mu\in G(X)}\left[\mobd(\mu;\sigma(G(0)))^{d_G}\right] \ \leq \ r(X) \quad\forall X\in\OM_n.
\]
For the moment, let us fix $X\in\OM_n$. For each $\mu\in\sigma(G(X))$, let $\lam_\mu$ be an 
eigenvalue of $G(0)$ such that
$\ilnhyper{\mu}{\lam_\mu}=\mobd(\mu;\sigma(G(0)))$. Now fix $\mu\in\sigma(G(X))$. The
above inequality leads to
\begin{equation}\label{E:hyperIneq1}
\hyper{\mu}{\lam_\mu} \ \leq \ r(X)^{1/d_G}.
\end{equation}
Applying Lemma~\ref{L:mobiusTrans} to the M{\"o}bius transformation
\[
T(z) \ = \ \frac{|\mu|z-\lam_\mu}{1-\overline{\lam_\mu}|\mu|z},
\]
we deduce that
\[
\hyper{\zt}{\lam_\mu} \ \geq \ \frac{||\mu|-|\lam_\mu||}{1-|\mu||\lam_\mu|}\quad\forall\zt:|\zt|=|\mu|.
\]
Applying the above fact to \eqref{E:hyperIneq1}, we get
\begin{align}
\frac{|\mu|-|\lam_\mu|}{1-|\mu||\lam_\mu|} \ &\leq \ r(X)^{1/d_G} \notag \\
\Rightarrow\quad |\mu| \ &\leq \frac{r(X)^{1/d_G}+|\lam_\mu|}{1+|\lam_\mu|r(X)^{1/d_G}},
\quad\mu\in\sigma(G(X)).
\label{E:specIneq}
\end{align}
Note that the function
\[
t\longmapsto\frac{r(X)^{1/d_G}+t}{1+r(X)^{1/d_G}t}, \quad t\geq 0,
\]
is an increasing function on $[0,\infty)$. Combining this fact with \eqref{E:specIneq}, we get
\[
|\mu| \ \leq \ \frac{r(X)^{1/d_G}+r(G(0))}{1+r(G(0))r(X)^{1/d_G}},
\]
which holds $\forall\mu\in\sigma(G(X))$, while the right-hand side is independent of $\mu$. Since
this is true for any arbitrary $X\in\OM_n$, we conclude that
\[
r(G(X)) \ \leq \ \frac{r(X)^{1/d_G}+r(G(0))}{1+r(G(0))r(X)^{1/d_G}} \quad\forall X\in\OM_n.
\]
\smallskip

In order to prove the sharpness of \eqref{E:SchwarzIneq}, let us fix an $n\geq 2$, and define
\[
\exep_n \ := \ \{A\in\OM_n:A \ \text{has a single eigenvalue of multiplicity $n$} \}.
\]
Pick any $d=1,\dots,n$, and define
\[
M_d(X) \ := \ \begin{cases}
                \quad [\tr(X)/n],        & \text{if $d=1$}, \\
                \ \begin{bmatrix} 
                        \ 0  & {} & {} & \tr(X)/n \ \\
                        \ 1  & 0  & {} & 0 \ \\
                        \ {} & \ddots & \ddots & \vdots \ \\
                        \ {} & {} & 1 & 0 \
                        \end{bmatrix}_{d\times d}, & \text{if $d\geq 2$},
                \end{cases}
\]
and, for the chosen $d$, define $\Shrp$ by the following block-diagonal matrix
\[
\Shrp(Y) \ := \begin{bmatrix}
                        \ M_d(X) & {} \ \\
                        \ {} & \dfrac{\tr(X)}{n}\mathbb{I}_{n-d} \
                        \end{bmatrix} \quad\forall X\in\OM_n.
\]
For our purposes $\shrp=\Shrp$ for each $A\in\exep_n$; i.e., the equality \eqref{E:sharp} will
will hold with the same function for each $A\in\exep_n$. To see this, note that
\begin{itemize}
\item $r(\Shrp(X))=|\tr(X)/n|^{1/d}$; and
\item $\Shrp(0)$ is nilpotent of degree $d$, whence $d_{\Shrp}=d$.
\end{itemize}
Therefore,
\[
\frac{r(A)^{1/d}+r(\Shrp(0))}{1+r(\Shrp(0))r(A)^{1/d}} \ = \
r(A)^{1/d} \ = \ r(\Shrp(A)) \quad\forall A\in\exep_n,
\]
which establishes \eqref{E:sharp}
\end{proof}
\end{custom}

\end{document}